# 复合材料板壳结构力学行为模拟和强度预测的统一高阶多尺度方法*

董灏[1)†] 葛步峰[2)] 高鸣远[3)]

1) (西安电子科技大学, 数学与统计学院, 西安 710071)

2) (西安电子科技大学, 电子工程学院, 西安 710071)

3) (西安电子科技大学, 人工智能学院, 西安 710071)

复合材料板壳结构复杂的细观结构导致直接模拟其力学问题需要巨大的计算量, 本文发展了可以有效模拟复合材料板壳结构力学行为并预测其强度的统一高阶多尺度方法. 首先, 通过对正交曲线坐标系下的多尺度弹性力学方程进行多尺度渐近分析, 建立了可以统一有效分析复合材料板和壳结构力学行为的高阶多尺度模型, 并得到了高阶多尺度解的误差估计. 然后, 结合材料强度理论, 建立复合材料板壳结构强度预测的高阶多尺度模型. 接下来, 基于建立的高阶多尺度模型, 建立了可以快速精确模拟复合材料板壳结构力学行为并预测其强度的多尺度算法. 最后通过数值实验验证了所建立的高阶多尺度方法的有效性, 数值实验结果表明高阶多尺度方法可以更加精确地捕捉复合材料板壳结构细观尺度的振荡行为. 本文建立的统一的高阶多尺度方法不仅适用于复合材料板壳结构力学性能的预测, 还可进一步推广应用于复合材料板壳结构多场耦合性能的预测.

**关键词**: 复合材料板壳结构, 多尺度渐近分析, 强度理论, 高阶多尺度方法, 误差分析

**PACS:** 72.80.Tm, 47.11.St, 46.50.+a, 02.60.-x

* 国家自然科学基金(批准号: 12001414)、中央高校基本科研业务费(批准号: JB210702)和新材料力学理论与应用湖北省重点实验室(武汉理工大学)开放课题(批准号: WUT-TAM202104)资助的课题.

† 通信作者. E-mail: donghao@mail.nwpu.edu.cn 电话: 13468772907

# 1 引 言

复合材料是由两种或多种不同性质的材料用物理或化学方法组成的新型材料, 具有优良的热学和力学性能, 并且集轻质、多功能及可设计性于一身, 具有广阔的工程应用前景[1]. 在实际工程应用中, 复合材料常被制造成多层板壳、夹层板壳、加筋板壳等板壳结构. 近年来随着航天、航空工业的快速发展, 复合材料板和壳结构由于其优异的力学性能被广泛应用于制造飞行器的机身、机翼、隔热结构等. 对于复合材料板壳结构, 从力学角度来看, 复合材料壳体的宏观结构整体尺寸远大于其细观结构特征尺寸, 直接进行数值求解需要足够细的计算网格才能捕捉到其细观尺度的信息. 从数学角度来看, 描述复合材料板壳结构力学行为的数学模型是一个具有强间断性和高频振荡性方程系数的偏微分方程系统, 直接数值求解计算量巨大. 为了给新型飞行器的结构设计和性能优化提供坚实的基础理论和高性能算法, 研究复合材料板壳结构力学行为模拟和强度预测的数值方法具有重要的应用价值和理论意义[2].

为了研究复合材料及其结构的物理、力学性能, 数学家和工程师发展了多种多尺度方法. 在连续介质力学范畴内, 主要存在着渐近均匀化方法(AHM)、异质多尺度方法(HMM)、变分多尺度方法(VMS)、多尺度有限元方法(MsFEM)、局部正交分解方法(LOD)、广义多尺度有限元方法(GMsFEM)、多尺度特征单元方法(MEM)、变分渐近均匀化方法(VAM)等[3–10]. 在这些多尺度方法中, 渐近均匀化方法具有严格的数学基础, 物理意义明确, 且能够很好地与有限元方法结合, 在实际工程应用中得到了广泛的应用. 针对传统均匀化方法和一阶渐近均匀化方法(一阶多尺度方法)无法有效捕捉复合材料局部振荡行为、数值精度低的缺点, 崔俊芝等针对周期型复合材料结构提出了高阶多尺度方法, 该方法能够同时保证对能量和动量平衡方程具有小尺度阶逼近性, 使得渐近均匀化方法从理论分析进入到数值计算和实际应用阶段[11–15]. 此外, 该高阶多尺度方法是一个扩展性很强的多尺度计算框架, 可以应用于连续介质力学范畴内的复合材料及其结构力学问题的模拟与分析. 对于复合材料板壳结构, Lewinski´ 等[16]和 Miara 等[17]使用渐近均匀化方法研究了复合材料壳体结构的物理、力学性能, 但他们所用到的细观单胞函数是定义在一个二维模型上的, 不能完全反映复合材料壳体结构的三维细观物理、力学行为, 更不能用于分析复合材料壳体结构细观尺度复杂的三维力学行为. Ma 等[18–19]基于渐近均匀化方法发展了轴对称复合材料圆柱壳和球对称复合材料球壳力学问题的二阶双尺度方法, 但建立的细观单胞模型仍是二维的, 无法用于分析计算一般复合材料壳体结构三维细观尺度的力学行为. 对于复合材料壳体结构力学行为的模拟问题, 力学界以往的研究多集中于发展层合壳的高阶壳理论以及低阶、高阶壳理论的应用[20–21], 所建立的壳体结构数学模型是一个二维模型, 提出的方法只能够计算复合材料壳体结构宏观尺度的力学响应, 无法分析复合材料壳体结构三维细观尺度的力学行为. 综上所述, 对适用于一般正交曲线坐标系统下复合材料板壳结构三维宏、细观力学行为模拟和强度预测的高阶多尺度方法的研究十分缺乏.

本文从正交曲线坐标系下描述复合材料板和壳结构力学行为的统一数学模型出发, 将复合材料板视为复合材料壳的特例, 建立了统一的复合材料板壳结构力学行为模拟和强度预测的高阶多尺度方法. 对于正交曲线坐标系下复合材料板壳结构多尺度力学问题, 从理论分析来看, 正交曲线坐标系下的几何方程不同于直角坐标系下的几何方程具有统一的形式, 这造成该问题双尺度分析的困难; 数值计算方面, 该问题的偏微分方程数学模型的强间断和快速振荡系数使得要精确捕捉复合材料板壳结构细观尺度振荡特性十分困难. 本文首先通过对正交曲线坐标系下多尺度弹性力学方程进行多尺度渐近分析, 建立了可以统一有效分析复合材料板和壳结构力学行为的高阶多尺度模型. 然后, 结合曲线坐标系下的泛函分析理论, 得到了高阶多尺度解的收敛性估计. 接下来, 结合材料强度理论, 建立了复合材料板壳结构强度预测的高阶多尺度模型. 最后,结合有限元方法和插值技术, 发展了可以快速精确模拟复合材料板壳结构力学行为并预测其强度的多尺度算法. 本文所建立的统一高阶多尺度方法具有很强的扩展性, 可以应用于多场耦合环境下复合材料板壳结构力学性能的计算和强度预测.

# 2 复合材料板壳结构力学行为模拟的高阶多尺度方法

本章首先建立复合材料板壳结构的统一数学模型, 然后对其进行多尺度渐近分析, 得到可以有效分析复合材料板壳结构的高阶多尺度模型. 然后, 结合曲线坐标系下的泛函分析理论, 对高阶多尺度模型的收敛性进行分析, 得到高阶多尺度解的误差估计.

## 2.1 复合材料板壳结构多尺度力学问题

基于板壳变形理论和多尺度理论[16,22_25], 建立描述复合材料板壳结构力学行为的正交曲线坐标下的偏微分方程数学模型如下所示:

$$\begin{cases} -\sum_{j=1}^{3}\frac{1}{H}\frac{\partial}{\partial\alpha_j}\left(\frac{H}{H_j}\sigma_{ij}^{\varepsilon}(\boldsymbol{\alpha})\right)-\sum_{j\neq i,j=1}^{3}\frac{1}{H_iH_j}\frac{\partial H_i}{\partial\alpha_j}\sigma_{ij}^{\varepsilon}(\boldsymbol{\alpha}) \\ \quad +\sum_{j\neq i,j=1}^{3}\frac{1}{H_iH_j}\frac{\partial H_j}{\partial\alpha_i}\sigma_{jj}^{\varepsilon}(\boldsymbol{\alpha})=f_i(\boldsymbol{\alpha}),\quad \text{in}\ \ \Omega, \\ \quad \boldsymbol{u}^{\boldsymbol{\varepsilon}}(\boldsymbol{\alpha})=\widehat{\boldsymbol{u}}(\boldsymbol{\alpha}),\quad \text{on}\ \ \partial\Omega_u, \\ \quad \sum_{j=1}^{3}\sigma_{ij}^{\varepsilon}(\boldsymbol{\alpha})\frac{n_j}{H_j}=\overline{\sigma_i}(\boldsymbol{\alpha}),\quad \text{on}\ \ \partial\Omega_\sigma. \end{cases} \tag{1}$$

其中$\boldsymbol{H}=\boldsymbol{H_1H_2H_3}$, $\boldsymbol{H_i}$表示$\boldsymbol{\alpha_i}$方向的拉梅系数; $\boldsymbol{\Omega}$是$\boldsymbol{R}^3$中的有界凸区域, $\boldsymbol{\partial\Omega}$为$\boldsymbol{\Omega}$的边界并且满足$\boldsymbol{\partial\Omega}=\boldsymbol{\partial\Omega_u}\cup\boldsymbol{\partial\Omega_\sigma}$; $\boldsymbol{\varepsilon}$代表复合材料板壳结构的小周期参数, $\boldsymbol{\alpha}=(\boldsymbol{\alpha_1},\boldsymbol{\alpha_2},\boldsymbol{\alpha_3})$为宏观尺度坐标, $\boldsymbol{f_i}(\boldsymbol{\alpha})$表示体力密度; $\widehat{\boldsymbol{u}}(\boldsymbol{\alpha})$ 和$\overline{\boldsymbol{\sigma}_{\boldsymbol{\iota}}}(\boldsymbol{\alpha})$分别表示边界$\boldsymbol{\partial\Omega_u}$和边界$\boldsymbol{\partial\Omega_\sigma}$上的位移和载荷分量.

复合材料板壳结构的平衡方程(1)中的本构方程为:

$$\sigma_{ij}^{\varepsilon}(\boldsymbol{\alpha}) = C_{ijkl}^{\varepsilon}(\boldsymbol{\alpha}) e_{kl}^{\varepsilon}(\boldsymbol{\alpha}). \tag{2}$$

几何方程为:

$$\begin{aligned}
e_{11}^{\varepsilon} &= \frac{1}{H_1}\frac{\partial u_1^{\varepsilon}}{\partial \alpha_1} + \frac{1}{H_1 H_2}\frac{\partial H_1}{\partial \alpha_2} u_2^{\varepsilon} + \frac{1}{H_1 H_3}\frac{\partial H_1}{\partial \alpha_3} u_3^{\varepsilon},\\
e_{22}^{\varepsilon} &= \frac{1}{H_2}\frac{\partial u_2^{\varepsilon}}{\partial \alpha_2} + \frac{1}{H_2 H_3}\frac{\partial H_2}{\partial \alpha_3} u_3^{\varepsilon} + \frac{1}{H_2 H_1}\frac{\partial H_2}{\partial \alpha_1} u_1^{\varepsilon},\\
e_{33}^{\varepsilon} &= \frac{1}{H_3}\frac{\partial u_3^{\varepsilon}}{\partial \alpha_3} + \frac{1}{H_3 H_1}\frac{\partial H_3}{\partial \alpha_1} u_1^{\varepsilon} + \frac{1}{H_3 H_2}\frac{\partial H_3}{\partial \alpha_2} u_2^{\varepsilon},\\
e_{12}^{\varepsilon} &= \frac{1}{2}\left[\frac{H_2}{H_1}\frac{\partial}{\partial \alpha_1}\left(\frac{u_2^{\varepsilon}}{H_2}\right) + \frac{H_1}{H_2}\frac{\partial}{\partial \alpha_2}\left(\frac{u_1^{\varepsilon}}{H_1}\right)\right],\\
e_{23}^{\varepsilon} &= \frac{1}{2}\left[\frac{H_3}{H_2}\frac{\partial}{\partial \alpha_2}\left(\frac{u_3^{\varepsilon}}{H_3}\right) + \frac{H_2}{H_3}\frac{\partial}{\partial \alpha_3}\left(\frac{u_2^{\varepsilon}}{H_2}\right)\right],\\
e_{13}^{\varepsilon} &= \frac{1}{2}\left[\frac{H_3}{H_1}\frac{\partial}{\partial \alpha_1}\left(\frac{u_3^{\varepsilon}}{H_3}\right) + \frac{H_1}{H_3}\frac{\partial}{\partial \alpha_3}\left(\frac{u_1^{\varepsilon}}{H_1}\right)\right].
\end{aligned} \tag{3}$$

通过求解复合材料板壳结构多尺度力学问题(1), 可以获得复合材料板壳结构的位移场$\boldsymbol{u_i^{\varepsilon}}(\boldsymbol{\alpha})$、应变场$\boldsymbol{e_{ij}^{\varepsilon}}(\boldsymbol{\alpha})$和应力场$\boldsymbol{\sigma_{ij}^{\varepsilon}}(\boldsymbol{\alpha})$, 对其力学行为进行模拟和分析.

## 2.2 力学行为模拟的高阶多尺度模型

定义复合材料板壳结构的细观尺度周期单胞$\boldsymbol{Q} = (\boldsymbol{0}, \boldsymbol{1})^{\boldsymbol{3}}$和细观尺度坐标$\boldsymbol{\beta} = \frac{\boldsymbol{\alpha}}{\boldsymbol{\varepsilon}} = \left(\frac{\boldsymbol{\alpha_1}}{\boldsymbol{\varepsilon}}, \frac{\boldsymbol{\alpha_2}}{\boldsymbol{\varepsilon}}, \frac{\boldsymbol{\alpha_3}}{\boldsymbol{\varepsilon}}\right) = (\boldsymbol{\beta_1}, \boldsymbol{\beta_2}, \boldsymbol{\beta_3})$, 根据复合材料板壳结构宏观尺度坐标和细观尺度坐标之间的尺度关系, 定义正交曲线坐标系下渐近均匀化方法的链式法则为

$$\frac{\partial}{\partial \alpha_i} = \frac{\partial}{\partial \alpha_i} + \frac{1}{\varepsilon}\frac{\partial}{\partial \beta_i}, \quad (i = 1,2,3). \tag{4}$$

对于复合材料板壳结构多尺度力学问题, 其几何方程缺乏直角坐标系下几何方程统一的形式, 我们首先假设复合材料板壳结构的应力具有如下多尺度渐近展开形式:

$$\sigma_{ij}^{\varepsilon}(\alpha) = \varepsilon^{-1}\sigma_{ij}^{(-1)}(\alpha, \beta) + \varepsilon^{0}\sigma_{ij}^{(0)}(\alpha, \beta) + \varepsilon^{1}\sigma_{ij}^{(1)}(\alpha, \beta) + O(\varepsilon^{2}). \tag{5}$$

然后将(5)代入(1)中并利用链式法则(2)对复合材料板壳结构弹性力学方程(1)进行多尺度渐近展开得

$$\begin{aligned}&\varepsilon^{-2}\left[-\sum_{j=1}^{3}\frac{1}{H}\frac{\partial}{\partial\alpha_j}\left(\frac{H}{H_j}\sigma_{ij}^{(-1)}\right)\right]+\varepsilon^{-1}\left[-\sum_{j=1}^{3}\frac{1}{H}\frac{\partial}{\partial\alpha_j}\left(\frac{H}{H_j}\sigma_{ij}^{(-1)}\right)\right.\\&\left.-\sum_{j=1}^{3}\frac{1}{H}\frac{\partial}{\partial\beta_j}\left(\frac{H}{H_j}\sigma_{ij}^{(0)}\right)-\sum_{j\neq i,j=1}^{3}\frac{1}{H_iH_j}\frac{\partial H_i}{\partial\alpha_j}\sigma_{ij}^{(-1)}+\sum_{j\neq i,j=1}^{3}\frac{1}{H_iH_j}\frac{\partial H_j}{\partial\alpha_i}\sigma_{jj}^{(-1)}\right]\\&+\varepsilon^{0}\left[-\sum_{j=1}^{3}\frac{1}{H}\frac{\partial}{\partial\alpha_j}\left(\frac{H}{H_j}\sigma_{ij}^{(0)}\right)-\sum_{j=1}^{3}\frac{1}{H}\frac{\partial}{\partial\beta_j}\left(\frac{H}{H_j}\sigma_{ij}^{(1)}\right)\right.\\&\left.-\sum_{j\neq i,j=1}^{3}\frac{1}{H_iH_j}\frac{\partial H_i}{\partial\alpha_j}\sigma_{ij}^{(0)}+\sum_{j\neq i,j=1}^{3}\frac{1}{H_iH_j}\frac{\partial H_j}{\partial\alpha_i}\sigma_{jj}^{(0)}\right]+\mathrm{O}(\varepsilon)=f_i.\end{aligned}\tag{6}$$

然后合并整理可得具有相同幂次小周期参数 ε 的方程如下:

$$O(\varepsilon^{-2}):-\sum_{j=1}^{3}\frac{1}{H}\frac{\partial}{\partial\beta_j}\left(\frac{H}{H_j}\sigma_{ij}^{(-1)}\right)=0.\tag{7}$$

$$\begin{aligned}O(\varepsilon^{-1}):&-\sum_{j=1}^{3}\frac{1}{H}\frac{\partial}{\partial\alpha_j}\left(\frac{H}{H_j}\sigma_{ij}^{(-1)}\right)-\sum_{j=1}^{3}\frac{1}{H}\frac{\partial}{\partial\beta_j}\left(\frac{H}{H_j}\sigma_{ij}^{(0)}\right)\\&-\sum_{j\neq i,j=1}^{3}\frac{1}{H_iH_j}\frac{\partial H_i}{\partial\alpha_j}\sigma_{ij}^{(-1)}+\sum_{j\neq i,j=1}^{3}\frac{1}{H_iH_j}\frac{\partial H_j}{\partial\alpha_i}\sigma_{jj}^{(-1)}=0.\end{aligned}\tag{8}$$

$$\begin{aligned}O(\varepsilon^{0}):&-\sum_{j=1}^{3}\frac{1}{H}\frac{\partial}{\partial\alpha_j}\left(\frac{H}{H_j}\sigma_{ij}^{(0)}\right)-\sum_{j=1}^{3}\frac{1}{H}\frac{\partial}{\partial\beta_j}\left(\frac{H}{H_j}\sigma_{ij}^{(1)}\right)\\&-\sum_{j\neq i,j=1}^{3}\frac{1}{H_iH_j}\frac{\partial H_i}{\partial\alpha_j}\sigma_{ij}^{(0)}+\sum_{j\neq i,j=1}^{3}\frac{1}{H_iH_j}\frac{\partial H_j}{\partial\alpha_i}\sigma_{jj}^{(0)}=f_i.\end{aligned}\tag{9}$$

为了简化下面的多尺度渐近分析过程, 定义新的宏观尺度和细观尺度微分算子如下所示:

$$\begin{cases}\psi_1=\dfrac{1}{H_1}\dfrac{\partial}{\partial\alpha_1},\psi_2=\dfrac{1}{H_2}\dfrac{\partial}{\partial\alpha_2},\psi_3=\dfrac{1}{H_3}\dfrac{\partial}{\partial\alpha_3},\\\widetilde{\psi_1}=\dfrac{1}{H_1}\dfrac{\partial}{\partial\beta_1},\widetilde{\psi_2}=\dfrac{1}{H_2}\dfrac{\partial}{\partial\beta_2},\widetilde{\psi_3}=\dfrac{1}{H_3}\dfrac{\partial}{\partial\beta_3}.\end{cases}\tag{10}$$

然后可建立如下新的多尺度链式法则

$$\psi_i=\psi_i+\varepsilon^{-1}\widetilde{\psi_i}\tag{11}$$

为了求出复合材料板壳结构的位移场, 我们假设位移场具有如下的多尺度渐近展开形式:

$$u_i^{\varepsilon}(\alpha)=u_i^{(0)}(\alpha,\beta)+\varepsilon u_i^{(1)}(\alpha,\beta)+\varepsilon^2u_i^{(2)}(\alpha,\beta)+\mathrm{O}(\varepsilon^3).\tag{12}$$

接下来, 我们对应变场进行如下多尺度渐近展开:

$$e_{ij}^{\varepsilon}(\alpha)=\varepsilon^{-1}e_{ij}^{(-1)}(\alpha,\beta)+\varepsilon^{0}e_{ij}^{(0)}(\alpha,\beta)+\varepsilon^{1}e_{ij}^{(1)}(\alpha,\beta)+O(\varepsilon^{2}). \tag{13}$$

结合(3)和(12), 可得应变场多尺度渐近展开项$e_{ij}^{(-1)}$、$e_{ij}^{(0)}$和$e_{ij}^{(1)}$的具体表达式为:

$$\begin{aligned}
e_{ij}^{(-1)}&=\frac{1}{2}[\tilde{\psi}_i(u_j^{(0)})+\tilde{\psi}_j(u_i^{(0)})],\\
e_{ij}^{(s)}&=e_{ij}^{(s)*}+\frac{1}{2}\left[\tilde{\psi}_i\left(u_j^{(s+1)}\right)+\tilde{\psi}_j\left(u_i^{(s+1)}\right)\right],\\
e_{11}^{(s)*}&=\psi_1\left(u_1^{(s)}\right)+\frac{\psi_2(H_1)}{H_1}u_2^{(s)}+\frac{\psi_3(H_1)}{H_1}u_3^{(s)},\\
e_{22}^{(s)*}&=\psi_2\left(u_2^{(s)}\right)+\frac{\psi_3(H_2)}{H_2}u_3^{(s)}+\frac{\psi_1(H_2)}{H_2}u_1^{(s)},\\
e_{33}^{(s)*}&=\psi_3\left(u_3^{(s)}\right)+\frac{\psi_1(H_3)}{H_3}u_1^{(s)}+\frac{\psi_2(H_3)}{H_3}u_2^{(s)},\\
e_{12}^{(s)*}&=\frac{1}{2}\left[\psi_1\left(u_2^{(s)}\right)+\psi_2\left(u_1^{(s)}\right)-\frac{\psi_1(H_2)}{H_2}u_2^{(s)}-\frac{\psi_2(H_1)}{H_1}u_1^{(s)}\right],\\
e_{23}^{(s)*}&=\frac{1}{2}\left[\psi_2\left(u_3^{(s)}\right)+\psi_3\left(u_2^{(s)}\right)-\frac{\psi_2(H_3)}{H_3}u_3^{(s)}-\frac{\psi_3(H_2)}{H_2}u_2^{(s)}\right],\\
e_{13}^{(s)*}&=\frac{1}{2}\left[\psi_1\left(u_3^{(s)}\right)+\psi_3\left(u_1^{(s)}\right)-\frac{\psi_1(H_3)}{H_3}u_3^{(s)}-\frac{\psi_3(H_1)}{H_1}u_1^{(s)}\right],\ s=0,1.
\end{aligned} \tag{14}$$

将(12)和(14)代入(5), 进一步可得$\sigma_{ij}^{(-1)}$、$\sigma_{ij}^{(0)}$和$\sigma_{ij}^{(1)}$的具体表达式为:

$$\begin{gathered}
\sigma_{ij}^{(-1)}=C_{ijkl}\widetilde{\psi_k}\left(u_l^{(0)}\right),\sigma_{ij}^{(0)}=C_{ijkl}e_{kl}^{(0)*}+C_{ijkl}\widetilde{\psi_k}\left(u_l^{(1)}\right),\\
\sigma_{ij}^{(1)}=C_{ijkl}e_{kl}^{(1)*}+C_{ijkl}\widetilde{\psi_k}\left(u_l^{(2)}\right).
\end{gathered} \tag{15}$$

下面根据(7) − (9), 依次求解位移场的各展开项的具体表达式. 首先, 将式(15)中的相关项代入式(7)可得:

$$\tilde{\psi}_j\left[C_{ijkl}\tilde{\psi}_k\left(u_l^{(0)}\right)\right]=0. \tag{16}$$

根据渐近均匀化理论可以得出$\boldsymbol{u}_{\boldsymbol{i}}^{(\mathbf{0})}$与细观尺度变量$\boldsymbol{\beta}$无关, 即

$$u_i^{(0)}(\alpha,\beta)=u_i^{(0)}(\alpha). \tag{17}$$

接下来, 将式(15)中的相关项代入式(8)可得:

$$\tilde{\psi}_j\left[C_{ijkl}\tilde{\psi}_k\left(u_l^{(1)}\right)\right]=-\tilde{\psi}_j(C_{ijkl})e_{kl}^{(0)*}. \tag{18}$$

根据式(18)中各变量之间的关系, 构造性地给出$\boldsymbol{u}_{\boldsymbol{i}}^{(\mathbf{1})}$的具体表达式为:

$$u_i^{(1)}(\boldsymbol{\alpha},\boldsymbol{\beta})=N_i^{mn}(H_1,H_2,H_3,\beta)e_{mn}^{(0)*},m,n=1,2,3 \tag{19}$$

式(19)中的$\boldsymbol{N_i^{mn}}(\boldsymbol{H_1},\boldsymbol{H_2},\boldsymbol{H_3},\boldsymbol{\beta})$表示一阶细观单胞函数, 将式(19)代入式(18)并附加周期边界条件可得$\boldsymbol{N_i^{mn}}(\boldsymbol{H_1},\boldsymbol{H_2},\boldsymbol{H_3},\boldsymbol{\beta})$满足如下的辅助单胞问题:

$$\begin{cases} \tilde{\psi}_j\left[C_{ijkl}\tilde{\psi}_k(N_l^{mn})\right] = -\tilde{\psi}_j\left(C_{ijmn}\right), \boldsymbol{\beta} \in Q, \\ N_i^{mn}(H_1,H_2,H_3,\beta)\text{关于细观变量}\beta\text{为}1-\text{周期函数}, \int_Q N_l^{mn} dQ = 0, \end{cases} \tag{20}$$

然后将式(15)代入式(9), 并对式(9)在细观单胞$\boldsymbol{Q}$上做均匀化积分平均, 可得宏观尺度均匀化力学问题如下:

$$\begin{cases} -\sum_{j=1}^{3}\frac{1}{H}\frac{\partial}{\partial \alpha_j}\left[\frac{H}{H_j}\left(\hat{C}_{ijkl}e_{kl}^{(0)*}\right)\right] - \sum_{j\neq i,j=1}^{3}\frac{1}{H_iH_j}\frac{\partial H_i}{\partial \alpha_j}\left[\hat{C}_{ijkl}e_{kl}^{(0)*}\right] \\ \quad + \sum_{j\neq i,j=1}^{3}\frac{1}{H_iH_j}\frac{\partial H_j}{\partial \alpha_i}\left[\hat{C}_{jjkl}e_{kl}^{(0)*}\right] = f_i, \quad \text{in } \Omega, \\ \boldsymbol{u}^{(0)}(\boldsymbol{\alpha}) = \hat{\boldsymbol{u}}(\boldsymbol{\alpha}), \quad \text{on } \partial\Omega_u, \\ \sum_{j=1}^{3}\hat{C}_{ijkl}\, e_{kl}^{(0)*}\frac{n_j}{H_j} = \bar{\sigma}_i(\alpha), \quad \text{on } \partial\Omega_\sigma. \end{cases} \tag{21}$$

其中宏观尺度均匀化材料参数$\widehat{C_{ijkl}}$的定义为:

$$\widehat{C_{ijkl}}(H_1,H_2,H_3) = \frac{1}{|Q|}\int_Q\left[C_{ijkl} + C_{ijmn}\tilde{\Psi}_m(N_n^{kl})\right]dQ\,. \tag{22}$$

复合材料板壳结构的宏观等效材料参数均依赖于宏观尺度变量$\boldsymbol{H_1}$、$\boldsymbol{H_2}$和$\boldsymbol{H_3}$, 多尺度方法将复合材料板壳结构在宏观尺度上转化成为一种特殊的功能梯度材料板壳结构.

下面用式(9)减去式(21)并将式(15)和式(19)代入可得高阶校正项$\boldsymbol{u_i^{(2)}}$满足如下的等式:

$$\widetilde{\Psi_j}\left[C_{ijkl}\tilde{\Psi}_k\left(u_l^{(2)}\right)\right]$$

$$= e_{mn}^{(0)*} * \begin{cases} -\sum_{j=1}^{3}\frac{H_j}{H}\psi_j\left[\frac{H}{H_j}\left(C_{ijmn} - \hat{C}_{ijmn} + C_{ijkl}\tilde{\psi}_k(N_l^{mn})\right)\right] \\ -\sum_{j\neq i,j=1}^{3}\left\{\frac{1}{H_i}\psi_j(H_i\}\right)\left(C_{ijmn} - \hat{C}_{ijmn} + C_{ijkl}\tilde{\psi}_k(N_l^{mn})\right) - \widetilde{\psi_j}\left(C_{ijkl}D_{klmn}\right) \\ + \sum_{j\neq i,j=1}^{3}\left\{\frac{1}{H_j}\psi_i(H_j\}\right)\left(C_{jjmn} - \hat{C}_{jjmn} + C_{jjkl}\tilde{\psi}_k(N_l^{mn})\right) \end{cases} \tag{22}$$

$$+\left[\hat{C}_{ijmn} - C_{ijmn} - C_{ijkl}\tilde{\Psi}_k(N_l^{mn}) - \tilde{\Psi}_l\left(C_{ilkj}N_k^{mn}\right)\right]\Psi_j\left(e_{mn}^{(0)*}\right).$$

其中$\boldsymbol{D_{klmn}}$的定义如下所示:

$$
\begin{aligned}
D_{11mn} &= \psi_1(N_1^{mn}) + \frac{\psi_2(H_1)}{H_1} N_2^{mn} + \frac{\psi_3(H_1)}{H_1} N_3^{mn}, \\
D_{22mn} &= \psi_2(N_2^{mn}) + \frac{\psi_3(H_2)}{H_2} N_3^{mn} + \frac{\psi_1(H_2)}{H_2} N_1^{mn}, \\
D_{33mn} &= \psi_3(N_3^{mn}) + \frac{\psi_1(H_3)}{H_3} N_1^{mn} + \frac{\psi_2(H_3)}{H_3} N_2^{mn}, \\
D_{12mn} &= \frac{1}{2}\left[\psi_1(N_2^{mn}) + \psi_2(N_1^{mn}) - \frac{\psi_1(H_2)}{H_2} N_2^{mn} - \frac{\psi_2(H_1)}{H_1} N_1^{mn}\right], \\
D_{23mn} &= \frac{1}{2}\left[\psi_2(N_3^{mn}) + \psi_3(N_2^{mn}) - \frac{\psi_2(H_3)}{H_3} N_3^{mn} - \frac{\psi_3(H_2)}{H_2} N_2^{mn}\right], \\
D_{13mn} &= \frac{1}{2}\left[\psi_1(N_3^{mn}) + \psi_3(N_1^{mn}) - \frac{\psi_1(H_3)}{H_3} N_3^{mn} - \frac{\psi_3(H_1)}{H_1} N_1^{mn}\right].
\end{aligned} \tag{24}
$$

根据式(23), 可以构造性地给出$\boldsymbol{u}_{\boldsymbol{i}}^{(\mathbf{2})}$的具体表达式为:

$$
u_i^{(2)}(\boldsymbol{\alpha},\boldsymbol{\beta}) = N_i^{jmn}(H_1,H_2,H_3,\beta)\psi_j\left(e_{mn}^{(0)*}\right) + W_i^{mn}(H_1,H_2,H_3,\beta)e_{mn}^{(0)*}. \tag{25}
$$

式(25)中的$N_i^{mn}(H_1,H_2,H_3,\beta)$和$W_i^{mn}(H_1,H_2,H_3,\beta)$表示二阶细观单胞函数, 将式(25)代入式(23)并附加周期边界条件可得$N_i^{mn}(H_1,H_2,H_3,\beta)$和$W_i^{mn}(H_1,H_2,H_3,\beta)$满足如下的辅助单胞问题:

$$
\begin{cases}
\tilde{\psi}_p\left[C_{ipkl}\tilde{\psi}_k\left(N_l^{jmn}\right)\right] = \hat{C}_{ijmn} - C_{ijmn} - C_{ijkl}\tilde{\psi}_k(N_l^{mn}) - \tilde{\psi}_l\left(C_{ilkj}N_k^{mn}\right), \beta \in Q, \\
N_l^{jmn}(H_1,H_2,H_3,\beta)\text{关于细观变量 β 为 1-周期函数}, \int_Q N_l^{jmn} dQ = 0.
\end{cases} \tag{26}
$$

$$
\begin{cases}
\tilde{\psi}_j\left[C_{ijkl}\tilde{\psi}_k(W_l^{mn})\right] = -\sum_{j=1}^{3}\frac{H_j}{H}\psi_j\left[\frac{H}{H_j}\left(C_{ijmn} - \hat{C}_{ijmn} + C_{ijkl}\tilde{\psi}_k(N_l^{mn})\right)\right] \\
-\sum_{j\neq i,j=1}^{3}\left\{\frac{1}{H_i}\psi_j(H_i\}\right)\left(C_{ijmn} - \hat{C}_{ijmn} + C_{ijkl}\tilde{\psi}_k(N_l^{mn})\right) - \widetilde{\psi_j}\left(C_{ijkl}D_{klmn}\right) \\
\quad + \sum_{j\neq i,j=1}^{3}\left\{\frac{1}{H_j}\psi_i\left(H_j\right\}\right)\left(C_{jjmn} - \hat{C}_{jjmn} + C_{jjkl}\tilde{\psi}_k(N_l^{mn})\right), \beta \in Q, \\
W_l^{mn}(H_1,H_2,H_3,\beta)\text{ 关于细观变量 β 为 1-周期函数}, \int_Q W_l^{mn} dQ = 0.
\end{cases} \tag{27}
$$

综上所述, 我们可以得到复合材料板壳结构多尺度力学问题(1)位移场的二阶双尺度解. 复合材料板壳结构的高阶多尺度模型包括一阶和二阶的细观尺度单胞函数观尺度的均匀化解和宏-细观耦合的高阶多尺度解. 通过建立高阶校正项, 该高阶多尺度模型可以统一有效分析复合材料板壳结构的力学行为.

*定理 1 复合材料板壳结构多尺度力学问题(1)的位移场具有如下的二阶双尺度解.*

$$
u_i^{\varepsilon}(\alpha) \approx u_i^{(0)} + \varepsilon N_i^{mn}(H_1,H_2,H_3,\beta)\varepsilon_{mn}^{(0)*} + \varepsilon^2\left[N_i^{jmn}(H_1,H_2,H_3,\beta)\psi_j\left(e_{mn}^{(0)*}\right) + W_i^{mn}(H_1,H_2,H_3,\beta)e_{mn}^{(0)*}\right]. \tag{28}
$$

**Remark 1** *令拉梅系数$H1 = 1, H2 = 1$和$H3 = 1$, 上述高阶多尺度模型退化为复合材料板的高阶多尺度模型; 令拉梅系数$H1 = 1, H2 = R2 + \alpha3$ 和$H3 = 1$, 其中$R2$为圆柱壳中面的沿$\alpha2$方向的曲率半径, 上述高阶多尺度模型退化为复合材料圆柱壳的高阶多尺度模型; 令拉梅系数$H1 = R1 + \alpha1, H2 = R2 + \alpha3$和$H3 = 1$, 其中$R1$和$R2$为双曲扁壳中面的沿$\alpha1$和$\alpha2$方向的曲率半径, 上述高阶多尺度模型退化为复合材料双曲扁壳的高阶多尺度模型.*

## 2.3 误差分析

分别定义复合材料板壳结构多尺度力学问题的一阶双尺度解$\boldsymbol{u}_{\boldsymbol{i}}^{(\mathbf{1\varepsilon})}$和二阶双尺度解$\boldsymbol{u}_{\boldsymbol{i}}^{(\mathbf{2\varepsilon})}$如下所示:

$$u_i^{(1\varepsilon)} = u_i^{(0)} + \varepsilon u_i^{(1)},\ \ u_i^{(2\varepsilon)} = u_i^{(0)} + \varepsilon u_i^{(1)} + \varepsilon^2 u_i^{(2)}. \tag{29}$$

然后定义一阶双尺度解的残差函数$\boldsymbol{u}_{\boldsymbol{\Delta i}}^{(\mathbf{1\varepsilon})}$和二阶双尺度解的残差函数$\boldsymbol{u}_{\boldsymbol{\Delta i}}^{(\mathbf{2\varepsilon})}$

$$u_{\Delta i}^{(1\varepsilon)} = u_i^{\varepsilon} - u_i^{(1\varepsilon)}, u_{\Delta i}^{(2\varepsilon)} = u_i^{\varepsilon} - u_i^{(2\varepsilon)}. \tag{30}$$

然后将残差函数$\boldsymbol{u}_{\boldsymbol{\Delta i}}^{(\mathbf{1\varepsilon})}$代入(1), 可得一阶双尺度解满足的残差方程:

$$\begin{cases} -\sum_{j=1}^{3}\frac{1}{H}\frac{\partial}{\partial\alpha_j}\left(\frac{H}{H_j}\sigma_{ij}^{\varepsilon}\left(\boldsymbol{u}_{\boldsymbol{\Delta}}^{(1\varepsilon)}\right)\right) - \sum_{j\neq i,j=1}^{3}\frac{1}{H_iH_j}\frac{\partial H_i}{\partial\alpha_j}\sigma_{ij}^{\varepsilon}\left(\boldsymbol{u}_{\boldsymbol{\Delta}}^{(1\varepsilon)}\right) \\ +\sum_{j\neq i,j=1}^{3}\frac{1}{H_iH_j}\frac{\partial H_j}{\partial\alpha_i}\sigma_{jj}^{\varepsilon}\left(\boldsymbol{u}_{\boldsymbol{\Delta}}^{(1\varepsilon)}\right) = F_{\mathbb{0}i}(\boldsymbol{\alpha},\boldsymbol{\beta}) + \varepsilon F_{\mathbb{1}i}(\boldsymbol{\alpha},\boldsymbol{\beta}), \text{in } \Omega, \\ u_{\Delta i}^{(1\varepsilon)}(\boldsymbol{\alpha}) = -\varepsilon N_i^{mn} e_{mn}^{(0)*} = \varepsilon\widehat{\chi_{1\iota}}(\boldsymbol{\alpha}),\ \text{on } \partial\Omega_u, \\ \sum_{j=1}^{3}\sigma_{ij}^{\varepsilon}\left(\boldsymbol{u}_{\boldsymbol{\Delta}}^{(1\varepsilon)}, T_{\Delta}^{(1\varepsilon)}\right)\frac{n_j}{H_j} = \sum_{j=1}^{3} I_{1ij}\frac{n_j}{H_j},\ \ \text{on } \partial\Omega_\sigma. \end{cases} \tag{31}$$

再将残差函数$\boldsymbol{u}_{\boldsymbol{\Delta i}}^{(\mathbf{2\varepsilon})}$代入(1), 可得二阶双尺度解满足的残差方程:

$$\begin{cases} -\sum_{j=1}^{3}\frac{1}{H}\frac{\partial}{\partial\alpha_j}\left(\frac{H}{H_j}\sigma_{ij}^{\varepsilon}\left(\boldsymbol{u}_{\boldsymbol{\Delta}}^{(2\varepsilon)}\right)\right) - \sum_{j\neq i,j=1}^{3}\frac{1}{H_iH_j}\frac{\partial H_i}{\partial\alpha_j}\sigma_{ij}^{\varepsilon}\left(\boldsymbol{u}_{\boldsymbol{\Delta}}^{(2\varepsilon)}\right) \\ +\sum_{j\neq i,j=1}^{3}\frac{1}{H_iH_j}\frac{\partial H_j}{\partial\alpha_i}\sigma_{jj}^{\varepsilon}\left(\boldsymbol{u}_{\boldsymbol{\Delta}}^{(2\varepsilon)}\right) = \varepsilon G_i(\boldsymbol{\alpha},\boldsymbol{\beta}),\ \ \text{in } \Omega, \\ u_{\Delta i}^{(2\varepsilon)}(\boldsymbol{\alpha}) = -\varepsilon N_i^{mn} e_{mn}^{(0)*} - \varepsilon^2 N_i^{jmn}\psi_j\left(e_{mn}^{(0)*}\right) - \varepsilon^2 W_i^{mn} e_{mn}^{(0)*} = \varepsilon\widehat{\chi_{2\iota}}(\boldsymbol{\alpha}),\ \text{on } \partial\Omega_u, \\ \sum_{j=1}^{3}\sigma_{ij}^{\varepsilon}\left(\boldsymbol{u}_{\boldsymbol{\Delta}}^{(2\varepsilon)}\right)\frac{n_j}{H_j} = \sum_{j=1}^{3} I_{2ij}\frac{n_j}{H_j},\ \ \text{on } \partial\Omega_\sigma. \end{cases} \tag{32}$$

从残差方程(31)和(32)可知, 一阶双尺度解对多尺度力学问题(1)的真解在点点意义下的逼近精度为$\boldsymbol{O}(\mathbf{1})$量级, 而二阶双尺度解对多尺度力学问题(1)真解的逼近精度为$\boldsymbol{O}(\boldsymbol{\varepsilon})$量级. 由于二阶双尺度解保证了力学平衡方程的小尺度逼近性, 二阶双尺度解能够精确捕捉复合材料细观尺度的振荡信息, 满足工程计算对数值精度的要求, 这正是发展高阶多尺度方法的最重要的原因.

下面给出二阶双尺度解在积分意义下的误差估计.

定理 2 设区域$\Omega$是具有 Lipschitz 连续边界的有界区域，$\boldsymbol{u}^{\varepsilon}(\boldsymbol{\alpha})$是多尺度力学问题(1)的解，则可得二阶双尺度解$\boldsymbol{u}^{(2\varepsilon)}(\boldsymbol{\alpha})$的误差估计：

$$||u^{\varepsilon}(\alpha)-u^{(2\varepsilon)}(\alpha)||_{\left(H^1(\Omega)\right)^3}\le C\varepsilon^{\frac{1}{2}}. \tag{33}$$

其中 C 为与小周期参数 ε 无关的常数.

证明: 首先定义类似于参考文献[26]中的截断函数$\boldsymbol{m}_{\boldsymbol{\varepsilon}}(\boldsymbol{\alpha})\in\boldsymbol{C}^{\infty}(\overline{\boldsymbol{\Omega}})$如下所示:

$$\begin{cases} m_{\varepsilon}(\boldsymbol{\alpha})=1,\ if\ dist\left(\boldsymbol{\alpha},\partial\Omega_u\right)\le\varepsilon,\\ m_{\varepsilon}(\boldsymbol{\alpha})=0,\ if\ dist\left(\boldsymbol{\alpha},\partial\Omega_u\right)\ge 2\varepsilon,\\ ||\nabla m_{\varepsilon}(\boldsymbol{\alpha})||_{L^{\infty}(\Omega)}\le C\varepsilon^{-1},\ if\ \varepsilon<dist\left(\boldsymbol{\alpha},\partial\Omega_u\right)<2\varepsilon. \end{cases} \tag{34}$$

然后定义新的残差函数$\boldsymbol{\delta}_{\boldsymbol{i}}^{(\mathbf{2}\boldsymbol{\varepsilon})}$如下所示:

$$\delta_i^{(2\varepsilon)}=u_{\Delta i}^{(2\varepsilon)}-\varepsilon m_{\varepsilon}(\alpha)\widehat{\chi_{2\iota}}(\alpha). \tag{35}$$

对残差方程(32)两端乘以$\boldsymbol{\delta}_{\boldsymbol{i}}^{(\mathbf{2}\boldsymbol{\varepsilon})}\boldsymbol{H}$再在$\boldsymbol{\Omega}$上做积分可得:

$$\int_{\Omega}\sum_{j=1}^{3}\frac{1}{H}\frac{\partial}{\partial\alpha_j}\left(\frac{H}{H_j}\sigma_{ij}^{\varepsilon}\left(\boldsymbol{u}_{\boldsymbol{\Delta}}^{(2\varepsilon)}\right)\right)\delta_i^{(2\varepsilon)}\,Hd\boldsymbol{\Omega}-\int_{\Omega}\sum_{j\neq i,j=1}^{3}\frac{1}{H_iH_j}\frac{\partial H_i}{\partial\alpha_j}\sigma_{ij}^{\varepsilon}\left(\boldsymbol{u}_{\boldsymbol{\Delta}}^{(2\varepsilon)}\right)\delta_i^{(2\varepsilon)}\,Hd\boldsymbol{\Omega}$$

$$+\int_{\Omega}\sum_{j\neq i,j=1}^{3}\frac{1}{H_iH_j}\frac{\partial H_j}{\partial\alpha_i}\sigma_{jj}^{\varepsilon}\left(\boldsymbol{u}_{\boldsymbol{\Delta}}^{(2\varepsilon)}\right)\delta_i^{(2\varepsilon)}H\,d\boldsymbol{\Omega}=\int_{\Omega}\varepsilon G_i(\boldsymbol{\alpha},\boldsymbol{\beta})\delta_i^{(2\varepsilon)}\,Hd\boldsymbol{\Omega}. \tag{36}$$

然后用$Green$公式化简式(36)可得:

$$-\int_{\Omega}\sigma_{ij}^{\varepsilon}\left(\boldsymbol{u}_{\boldsymbol{\Delta}}^{(2\varepsilon)}\right)\psi_j\left(\delta_i^{(2\varepsilon)}\right)H\,d\boldsymbol{\Omega}-\int_{\Omega}\sum_{j\neq i,j=1}^{3}\frac{1}{H_iH_j}\frac{\partial H_i}{\partial\alpha_j}\sigma_{ij}^{\varepsilon}\left(\boldsymbol{u}_{\boldsymbol{\Delta}}^{(2\varepsilon)}\right)\delta_i^{(2\varepsilon)}\,Hd\boldsymbol{\Omega} \tag{37}$$

$$+\int_{\Omega}\sum_{j\neq i,j=1}^{3}\frac{1}{H_iH_j}\frac{\partial H_j}{\partial\alpha_i}\sigma_{jj}^{\varepsilon}\left(\boldsymbol{u}_{\boldsymbol{\Delta}}^{(2\varepsilon)}\right)\delta_i^{(2\varepsilon)}\,Hd\boldsymbol{\Omega}+\int_{\partial\Omega_\sigma}\sum_{j=1}^{3}\frac{1}{H_j}I_{2ij}\,n_j\delta_i^{(2\varepsilon)}H\,ds=\int_{\Omega}\varepsilon G_i(\boldsymbol{\alpha},\boldsymbol{\beta})\delta_i^{(2\varepsilon)}\,Hd\boldsymbol{\Omega}$$

再根据正交曲线坐标系下应变的定义(3), 式(37)可进一步简化为

$$-\int_{\Omega}C_{ijkl}^{\varepsilon}e_{kl}^{\varepsilon}\left(\boldsymbol{u}_{\boldsymbol{\Delta}}^{(2\varepsilon)}\right)e_{ij}^{\varepsilon}\left(\boldsymbol{\delta}^{(2\varepsilon)}\right)H\,d\boldsymbol{\Omega}+\int_{\partial\Omega_\sigma}\sum_{j=1}^{3}\frac{1}{H_j}I_{2ij}\,n_j\delta_i^{(2\varepsilon)}\,Hds=\int_{\Omega}\varepsilon G_i(\boldsymbol{\alpha},\boldsymbol{\beta})\delta_i^{(2\varepsilon)}\,Hd\boldsymbol{\Omega} \tag{38}$$

然后, 我们替换$\boldsymbol{u}_{\boldsymbol{\Delta}}^{(2\varepsilon)}$为$\boldsymbol{\delta}^{(2\varepsilon)}+\boldsymbol{\varepsilon m_{\varepsilon}(\alpha)\widehat{\chi_2}(\alpha)}$可以得到:

$$\int_{\Omega}C_{ijkl}^{\varepsilon}e_{kl}^{\varepsilon}\left(\boldsymbol{\delta}^{(2\varepsilon)}\right)e_{ij}^{\varepsilon}\left(\boldsymbol{\delta}^{(2\varepsilon)}\right)Hd\boldsymbol{\Omega}=\int_{\partial\Omega_\sigma}\sum_{j=1}^{3}\frac{1}{H_j}I_{2ij}\,n_j\delta_i^{(2\varepsilon)}H\,ds \tag{39}$$

$$-\int_{\Omega}C_{ijkl}^{\varepsilon}e_{kl}^{\varepsilon}\left(\varepsilon\,m_{\varepsilon}(\boldsymbol{\alpha})\widehat{\boldsymbol{\chi}_2}(\boldsymbol{\alpha})\right)e_{ij}^{\varepsilon}\left(\boldsymbol{\delta}^{(2\varepsilon)}\right)\mathrm{H}d\boldsymbol{\Omega}-\int_{\Omega}\varepsilon\,G_i(\boldsymbol{\alpha},\boldsymbol{\beta})\delta_i^{(2\varepsilon)}\,\mathrm{H}d\boldsymbol{\Omega}$$

下面通过对式(39)两端分别进行估计, 最终得到误差估计(33). 首先, 利用曲线坐标系的$Korn$不等式[24]对式(39)左端进行估计可得如下不等式:

$$\left|\int_{\Omega} C_{ijkl}^{\varepsilon} e_{kl}^{\varepsilon}(\boldsymbol{\delta}^{(2\varepsilon)}) e_{ij}^{\varepsilon}(\boldsymbol{\delta}^{(2\varepsilon)}) H d\boldsymbol{\Omega}\right| \geq C|\delta^{(2\varepsilon)}|_{(H^1(\Omega))^3}^2 \tag{40}$$

接下来, 根据$H\ddot{o}lder$不等式、均值不等式和参考文献 1.2, 对式(39)右端进行估计可得如下不等式:

$$|\int_{\partial\Omega_\sigma} \sum_{j=1}^{3} \frac{1}{H_j} I_{2ij} n_j \delta_i^{(2\varepsilon)} H\, ds - \int_{\Omega} C_{ijkl}^{\varepsilon} e_{kl}^{\varepsilon}(\varepsilon m_\varepsilon(\boldsymbol{\alpha})\widehat{\chi_2}(\alpha))\, e_{ij}^{\varepsilon}(\delta^{(2\varepsilon)}) H d\boldsymbol{\Omega} - \int_{\Omega} \varepsilon G_i(\alpha,\beta)\delta_i^{(2\varepsilon)} H d\boldsymbol{\Omega}|$$

$$\leq C\varepsilon^{\frac{1}{2}}||\delta^{(2\varepsilon)}||_{(H^1(\Omega))^3} + C||\varepsilon m_\varepsilon(\alpha)\widehat{\chi_2}(\alpha)||_{(H^1(\Omega))^3}\left|\left|\delta^{(2\varepsilon)}\right|\right|_{(H^1(\Omega))^3} + C\left|\left|\varepsilon G_i(\alpha,\beta)\right|\right|_{L^2(\Omega)}||\delta_i^{(2\varepsilon)}||_{L^2(\Omega)} \tag{41}$$

$$\leq C\varepsilon^{\frac{1}{2}}||\boldsymbol{\delta}^{(2\varepsilon)}||_{(H^1(\Omega))^3} + C\varepsilon||m_\varepsilon(\boldsymbol{\alpha})\widehat{\boldsymbol{\chi}_2}(\boldsymbol{\alpha})||_{(H^1(\Omega))^3}||\boldsymbol{\delta}^{(2\varepsilon)}||_{(H^1(\Omega))^3} + C\varepsilon||\boldsymbol{\delta}^{(2\varepsilon)}||_{(H^1(\Omega))^3}$$

结合式(40)和式(41)可得如下不等式:

$$||\boldsymbol{\delta}^{(2\varepsilon)}||_{(H^1(\Omega))^3} \leq C\varepsilon^{\frac{1}{2}} + C\varepsilon||m_\varepsilon(\boldsymbol{\alpha})\widehat{\boldsymbol{\chi}_2}(\boldsymbol{\alpha})||_{(H^1(\Omega))^3} + C\varepsilon. \tag{42}$$

再利用式(35)并对(42)左端使用三角不等式可得:

$$\left||u_\Delta^{(2\varepsilon)}\right||_{(H^1(\Omega))^3} \leq C\varepsilon^{\frac{1}{2}} + 2C\varepsilon\left||m_\varepsilon(\alpha)\widehat{\chi_2}(\alpha)\right||_{(H^1(\Omega))^3} + C\varepsilon \tag{43}$$

然后利用参考文献[26]中对截断函数的误差估计, 可得如下的不等式:

$$||m_\varepsilon(\alpha)\widehat{\chi_2}(\alpha)||_{(H^1(\Omega))^3} = ||m_\varepsilon(\alpha)\widehat{\chi_2}(\alpha)||_{(H^1(K_\varepsilon))^3} \leq C\varepsilon^{-\frac{1}{2}} \tag{44}$$

其中$\boldsymbol{K_\varepsilon} = \{\boldsymbol{\alpha}|\boldsymbol{dist}(\boldsymbol{\alpha}, \boldsymbol{\partial\Omega_u}) \leq \boldsymbol{2\varepsilon}\} \cap \boldsymbol{\Omega}$. 最后, 结合式(43)和式(44)可得最终误差估计:

$$\left||\boldsymbol{u}^\varepsilon(\boldsymbol{\alpha}) - \boldsymbol{u}^{(2\varepsilon)}(\boldsymbol{\alpha})\right||_{(H^1(\Omega))^3} = \left||\boldsymbol{u}_{\boldsymbol{\Delta}}^{(2\varepsilon)}\right||_{(H^1(\Omega))^3} \leq C\varepsilon^{\frac{1}{2}}. \tag{45}$$

# 3 复合材料板壳强度预测的高阶多尺度模型

复合材料板壳结构的强度是其众多力学性能中非常重要的一个参数, 准确预测其强度对研究复合材料板壳结构在复杂应力条件下屈服或产生破坏的规律具有十分重要的意义. 本节首先建立复合材料板壳结构应变和应力场的高阶多尺度计算公式, 然后基于得到的高精度应变和应力场计算公式, 结合材料强度理论, 建立复合材料板壳结构的强度预测公式.

## 3.1 复合材料板壳结构应变和应力场的高阶多尺度计算公式

基于复合材料板壳结构多尺度力学问题的高阶多尺度模型(28), 我们可以得到复合材料板壳结构应变场的高精度高阶多尺度计算公式:

$$
\begin{aligned}
\widehat{e_{11}^{\varepsilon}} &= \frac{1}{H_1}\frac{\partial u_1^{(2\varepsilon)}}{\partial \alpha_1} + \frac{1}{H_1 H_2}\frac{\partial H_1}{\partial \alpha_2} u_2^{(2\varepsilon)} + \frac{1}{H_1 H_3}\frac{\partial H_1}{\partial \alpha_3} u_3^{(2\varepsilon)}, \\
\widehat{e_{22}^{\varepsilon}} &= \frac{1}{H_2}\frac{\partial u_2^{(2\varepsilon)}}{\partial \alpha_2} + \frac{1}{H_2 H_3}\frac{\partial H_2}{\partial \alpha_3} u_3^{(2\varepsilon)} + \frac{1}{H_2 H_1}\frac{\partial H_2}{\partial \alpha_1} u_1^{(2\varepsilon)}, \\
\widehat{e_{33}^{\varepsilon}} &= \frac{1}{H_3}\frac{\partial u_3^{(2\varepsilon)}}{\partial \alpha_3} + \frac{1}{H_3 H_1}\frac{\partial H_3}{\partial \alpha_1} u_1^{(2\varepsilon)} + \frac{1}{H_3 H_2}\frac{\partial H_3}{\partial \alpha_2} u_2^{(2\varepsilon)}, \\
\widehat{e_{12}^{\varepsilon}} &= \frac{1}{2}\left[\frac{H_2}{H_1}\frac{\partial}{\partial \alpha_1}\left(\frac{u_2^{(2\varepsilon)}}{H_2}\right) + \frac{H_1}{H_2}\frac{\partial}{\partial \alpha_2}\left(\frac{u_1^{(2\varepsilon)}}{H_1}\right)\right] \\
\widehat{e_{23}^{\varepsilon}} &= \frac{1}{2}\left[\frac{H_3}{H_2}\frac{\partial}{\partial \alpha_2}\left(\frac{u_3^{(2\varepsilon)}}{H_3}\right) + \frac{H_2}{H_3}\frac{\partial}{\partial \alpha_3}\left(\frac{u_2^{(2\varepsilon)}}{H_2}\right)\right], \\
\widehat{e_{13}^{\varepsilon}} &= \frac{1}{2}\left[\frac{H_3}{H_1}\frac{\partial}{\partial \alpha_1}\left(\frac{u_3^{(2\varepsilon)}}{H_3}\right) + \frac{H_1}{H_3}\frac{\partial}{\partial \alpha_3}\left(\frac{u_1^{(2\varepsilon)}}{H_1}\right)\right].
\end{aligned}
\tag{46}
$$

进一步, 我们可以得到复合材料板壳结构应力场的高精度高阶多尺度计算公式:

$$
\widehat{\sigma_{ij}^{\varepsilon}}(\alpha) = C_{ijkl}^{\varepsilon}(\alpha)\widehat{e_{kl}^{\varepsilon}}(\alpha). \tag{47}
$$

## 3.2 复合材料板壳结构的强度预测公式

复合材料板壳结构是由复合材料生产制备而成, 它的强度依赖于其组分材料的屈服状态, 对复合材料板壳结构的强度预测相比均质板壳结构更加复杂. 在过去的几十年里, 学者们相继发展了多种不同的强度准则, 每一种强度准则(或称屈服准则)通常只适用于某一类材料, 但对于复合材料板壳结构的强度预测没有有效的预测准则.

为了有效地预测复合材料板壳结构的强度, 对于复合材料板壳结构的不同组分材料根据其材料性质选择合适的强度准则. 然后, 给定一个载荷, 利用复合材料板壳结构应变和应力场的高阶多尺度计算公式计算复合材料板壳结构的应变和应力场, 并根据各组分材料的选用的强度准则, 判断复合材料板壳结构各组分材料的屈服状态. 如果没有满足强度准则, 继续增加载荷, 直到找到复合

材料板壳结构的临界载荷$\boldsymbol{S}(\Omega)$，即在此载荷下, 复合材料板壳结构内应力(应变)最大的单元正好处于破坏临界点上. 本文采用米塞斯($Von-Mises$)等效应力强度准则[27-28], 形式如下所示:

$$\sigma_e=\frac{1}{\sqrt{2}}\sqrt{(\sigma_{11}-\sigma_{22})^2+(\sigma_{22}-\sigma_{33})^2+(\sigma_{33}-\sigma_{11})^2+6(\sigma_{12}^2+\sigma_{23}^2+\sigma_{31}^2)}<S_e \quad (48)$$

其中$\boldsymbol{\sigma_e}$为米塞斯($Von-Mises$)等效应力, $\boldsymbol{S_e}$为材料通过单向拉伸(压缩)试验测得的弹性强度极限值. 通过将复合材料板壳结构的应力计算结果代入(48)进行等式替换$\boldsymbol{\sigma_{ij}}=\widehat{\boldsymbol{\sigma}_{ij}^{\varepsilon}}(\boldsymbol{\alpha})$, 当且仅当上述不等式不成立时,复合材料板壳结构中对应的组分材料处于屈服状态, 达到弹性极限强度, 复合材料板壳结构开始发生破坏。

# 4 力学行为模拟和强度预测的多尺度算法

基于复合材料板壳结构力学行为模拟和强度预测的高阶多尺度模型, 给出相应的多尺度算法如下:

(S1) 确定细观单胞$\boldsymbol{Q}$的几何构造以及复合材料中各相材料的材料参数, 定义$\boldsymbol{J^{h_1}}=\{\boldsymbol{K}\}$为细观单胞$\boldsymbol{Q}$的一族正则化四面体网格剖分, 再定义$\boldsymbol{V_{h_1}}(\boldsymbol{Q})$为细观单胞$\boldsymbol{Q}$上的协调有限元空间.

(S2) 在宏观区域$\Omega$上选取一定数量的代表点$\boldsymbol{\alpha_I}$ ,然后根据细观单胞问题(20)求解代表点$\boldsymbol{\alpha_I}$处的依赖于宏观尺度参数$\big(H_1(\boldsymbol{\alpha_I}),H_2(\boldsymbol{\alpha_I}),H_3(\boldsymbol{\alpha_I})\big)$的一阶细观单胞函数, 其有限元求解方案如下:

$$N_i^{mn}(H_1,H_2,H_3,\beta)-\int_Q C_{ijkl}\widetilde{\Psi}_k(N_l^{mn})\widetilde{\Psi}_j\big(\upsilon_i^{h_1}\big)dQ=\int_Q C_{ijmn}\widetilde{\Psi}_j\big(\upsilon_i^{h_1}\big)dQ\,,\ \forall\upsilon^{h_1}\in\Big(V_{h_1}(Q)\Big)^3. \quad (49)$$

(S3) 确定宏观均匀化区域$\Omega$的几何构造, 定义$\boldsymbol{J^{h_0}}=\{\boldsymbol{e}\}$为宏观均匀化区域$\Omega$的一族正则化四面体网格剖分,再定义$\boldsymbol{V_{h_0}}(\boldsymbol{\Omega})$为宏观均匀化区域$\Omega$上的协调有限元空间.

(S4) 通过插值方法, 求解得到有限元空间$\boldsymbol{V_{h_0}}(\boldsymbol{\Omega})$上每个具有宏观尺度参数$H_1(\boldsymbol{\alpha_I}),H_2(\boldsymbol{\alpha_I}),H_3(\boldsymbol{\alpha_I})$的单元节点$\alpha$处的宏观尺度均匀化材料参数$\widehat{C_{ijkl}}(H_1,H_2,H_3)$. 然后给出宏观均匀化力学问题(21)的有限元求解方案如下:

$$\begin{cases}\iint_\Omega \widehat{C_{ijkl}}e_{kl}^{(0)*}\big(\boldsymbol{u^{(0)}}\big)e_{ij}^{(0)*}\big(\boldsymbol{v^{h_0}}\big)Hd\boldsymbol{\Omega}=\int_\Omega f_i v_i^{h_0}\,Hd\boldsymbol{\Omega}+\int_{\partial\Omega_\sigma}\overline{\sigma_i}(\boldsymbol{\alpha})v_i^{h_0}\,Hds,\ \forall\boldsymbol{v^{h_0}}\in\Big(V_{h_0}(\boldsymbol{\Omega})\Big)^3, \\ \boldsymbol{u^{(0)}}(\boldsymbol{\alpha})=\widehat{\boldsymbol{u}}(\boldsymbol{\alpha}),\quad \text{on }\partial\Omega_u.\end{cases} \quad (50)$$

(S5) 使用求解一阶单胞函数的计算网格$\boldsymbol{J^{h_1}}=\{\boldsymbol{K}\}$和有限元空间$\boldsymbol{V_{h_1}}(\boldsymbol{Q})$, 根据细观单胞问题(25)和(26)求解代表点$\boldsymbol{\alpha_I}$处的依赖于宏观尺度参数$\big(\boldsymbol{H_1}(\boldsymbol{\alpha_I}),\boldsymbol{H_2}(\boldsymbol{\alpha_I}),\boldsymbol{H_3}(\boldsymbol{\alpha_I})\big)$的二阶细观单胞函数$\boldsymbol{N_i^{jmn}}\big(\boldsymbol{H_1},\boldsymbol{H_2},\boldsymbol{H}_{3},,\boldsymbol{\beta}\big)$以及$\boldsymbol{W_i^{jmn}}\big(\boldsymbol{H_1},\boldsymbol{H_2},\boldsymbol{H}_{3},,\boldsymbol{\beta}\big)$, 二阶细观单胞函数的有限元计算方案类似于一阶细观单胞函数的有限元计算方案, 不再赘述.

(S6) 对于任意点$\boldsymbol{\alpha}=(\boldsymbol{\alpha_1},\boldsymbol{\alpha_2},\boldsymbol{\alpha_3})\in\Omega$, 使用插值方法求解得到该点处的一阶、二阶细观单胞函数和宏观均匀化解, 再利用单元平均方法[29]计算得到$\boldsymbol{\psi_j}\left(\boldsymbol{e}_{\boldsymbol{mn}}^{(\boldsymbol{0})*}\right)$和$\boldsymbol{e}_{\boldsymbol{mn}}^{(\boldsymbol{0})*}$, 然后根据(27)计算得到任意点$\alpha$处多尺度力学问题(1)位移场的二阶双尺度有限元解. 进一步, 可以根据(46)和(47)计算得到应变和应力场的二阶双尺度有限元解. 此外, 我们还可以应用后处理技术进一步提高复合材料板壳结构的位移、应变和应力的数值精度[30].

(S7) 采用迭代思想, 逐渐增加施加的荷载, 根据第 3.2 节给出的强度预测公式对复合材料板壳结构进行评估计算出复合材料板壳结构的屈服强度. 此外应该强调的是, 我们采用二分法来搜索和确定合适的用于强度预测的初始载荷, 可以提高复合材料板壳结构强度预测的计算效率.

Remark 2 本文中复合材料板壳结构的几何建模和多尺度计算均通过 Freefem++软件二次开发实现. Freefem++是一款免费、开放源代码的有限元计算软件, 由巴黎第六大学研究人员开发, 具有良好的有限元网格剖分和二次开发功能.

# 5 数值实验

本节通过数值实验来验证所发展复合材料板壳结构力学行为模拟和强度预测的统一高阶多尺度方法的有效性.

## 5.1 复合材料板壳结构力学行为模拟

本节对由正交编织复合材料、曲线编织复合材料和四向编织复合材料组成的复合材料板和壳的力学行为进行模拟, 研究的复合材料板壳结构的宏、细观结构如图 1 所示. 复合材料板的宏观结构$\boldsymbol{\Omega}=(\boldsymbol{\alpha_1},\boldsymbol{\alpha_2},\boldsymbol{\alpha_3})=[\boldsymbol{0},\boldsymbol{1}]\times[\boldsymbol{0},\boldsymbol{1}]\times[\boldsymbol{0},\boldsymbol{0.2}]\boldsymbol{cm^3}$且有小周期参数$\boldsymbol{\varepsilon}=\boldsymbol{1/25}$. 复合材料壳的宏观结构$\boldsymbol{\Omega}=(\boldsymbol{\alpha_1},\boldsymbol{\alpha_2},\boldsymbol{\alpha_3})=[-\boldsymbol{\pi/9},\boldsymbol{\pi/9}]\times[-\boldsymbol{\pi/9},\boldsymbol{\pi/9}]\times[-\boldsymbol{\pi/54},\boldsymbol{\pi/54}]\boldsymbol{cm}$, 且有小周期参数$\boldsymbol{\varepsilon}=\boldsymbol{1/128}$.

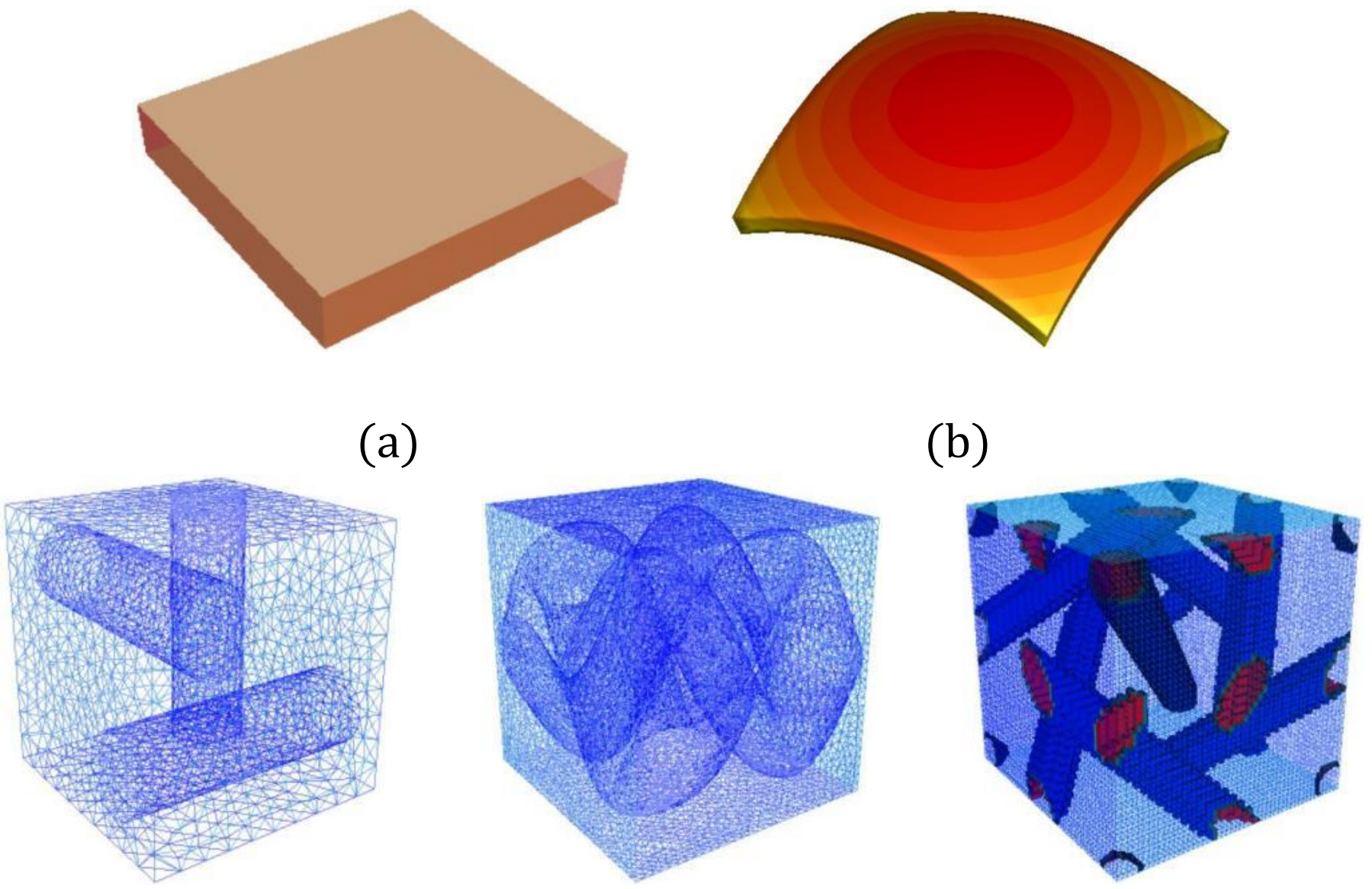

(a) (b)

(c) (d) (e)

图 1 复合材料板和壳宏、细观尺度结构示意图: (a) 复合材料板宏观结构; (b) 复合材料壳宏观结构; (c) 细观单胞 I; (d) 细观单胞 II; (e) 细观单胞 III.

Fig. 1. The macroscopic and mesoscopic configurations of composite plate and shell: (a) macroscopic structure of composite plate; (b) macroscopic structure of composite shell; (c) mesoscopic unit cell I; (d) mesoscopic unit cell II; (e) mesoscopic unit cell III.

我们研究的复合材料板壳结构的基体杨氏模量为 11700.00GPa, 泊松比为 0.321. 夹杂的杨氏模量为 6.62GPa, 泊松比为 0.333. 对复合材料板壳结构垂直于$\alpha_1$和$\alpha_2$方向的四个面固支, 并施加体力$(f_1, f_2, f_3) = (0,0,-10000.0)N/cm^3$. 然后对复合材料板壳宏观均匀化结构和细观单胞分别进行四面体网格剖分, 表 1 为各套网格的单元和节点信息.

表 1 计算资源信息.

Table 1. Information of computational cost.

| 有限元计算网格 | 复合材料板均匀化结构 | 复合材料壳均匀化结构 | 细观单胞 I | 细观单胞 II | 细观单胞 III |
|---|---|---|---|---|---|
| 单元个数 | 1200000 | 884736 | 22030 | 111551 | 146334 |
| 节点个数 | 214221 | 159953 | 4202 | 20759 | 27000 |

多尺度模拟完成后, 得到的数值结果如图 2 和 3 所示.

从数值模拟结果图 2 和 3 可以发现, 二阶双尺度解相比均匀化解与一阶双尺度解能够精确地捕捉到复合

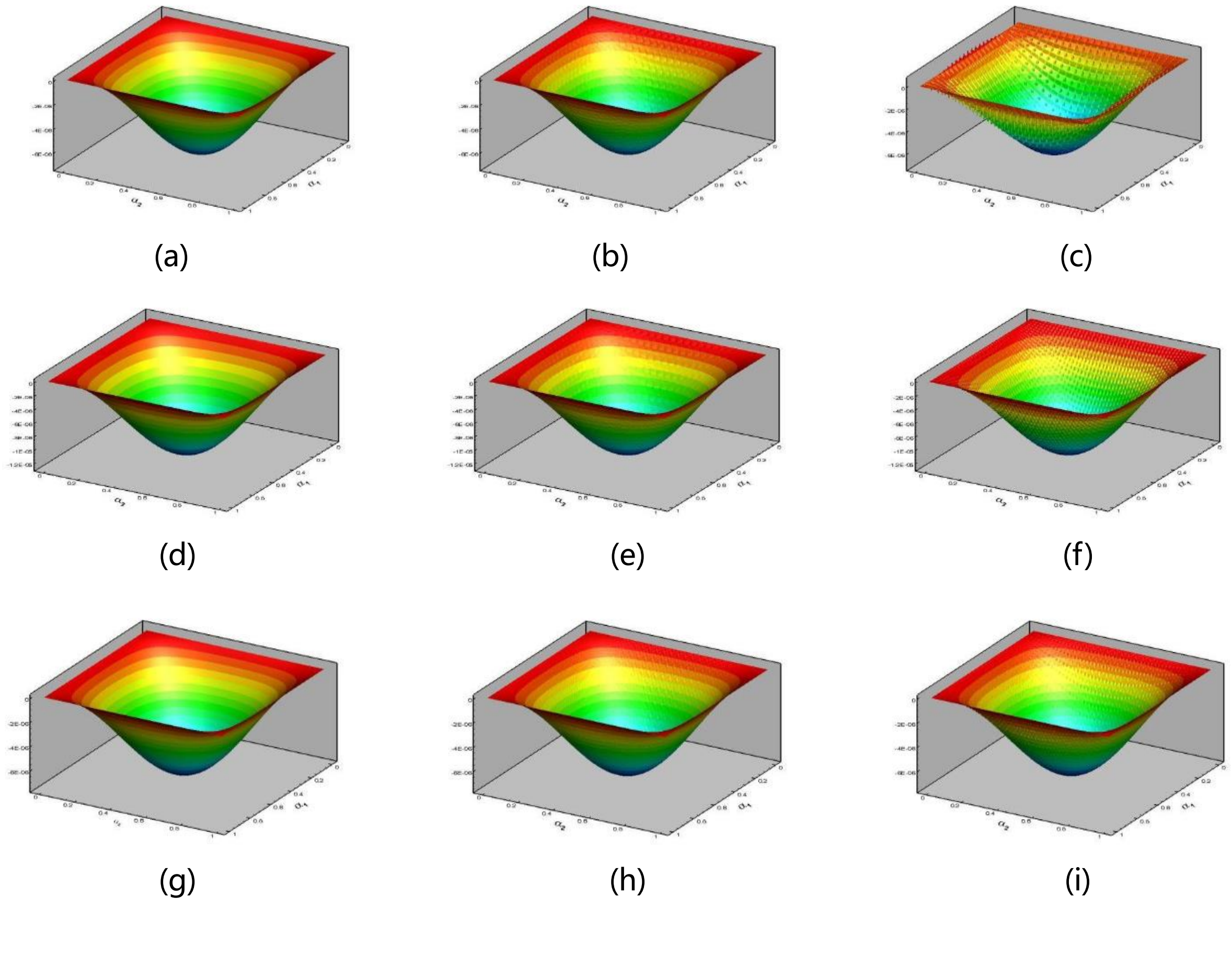

图 2 复合材料板位移场计算结果示意图, 细观单胞 I: (a) $\boldsymbol{u}_3^{(0)}$; (b) $\boldsymbol{u}_3^{(1\varepsilon)}$; (c) $\boldsymbol{u}_3^{(2\varepsilon)}$; 细观单胞 II: (d)$\boldsymbol{u}_3^{(0)}$; (e) $\boldsymbol{u}_3^{(1\varepsilon)}$; (f) $\boldsymbol{u}_3^{(2\varepsilon)}$; 细观单胞 III: (g) $\boldsymbol{u}_3^{(0)}$; (h) $\boldsymbol{u}_3^{(1\varepsilon)}$; (i) $\boldsymbol{u}_3^{(2\varepsilon)}$.

Fig. 2. The computational results of displacement field of composite plate, mesoscopic configuration I: (a) $\boldsymbol{u}_3^{(0)}$; (b) $\boldsymbol{u}_3^{(1\varepsilon)}$; (c) $\boldsymbol{u}_3^{(2\varepsilon)}$; mesoscopic configuration II: (d)$\boldsymbol{u}_3^{(0)}$; (e) $\boldsymbol{u}_3^{(1\varepsilon)}$; (f) $\boldsymbol{u}_3^{(2\varepsilon)}$; mesoscopic configuration III: (g) $\boldsymbol{u}_3^{(0)}$; (h) $\boldsymbol{u}_3^{(1\varepsilon)}$; (i) $\boldsymbol{u}_3^{(2\varepsilon)}$.

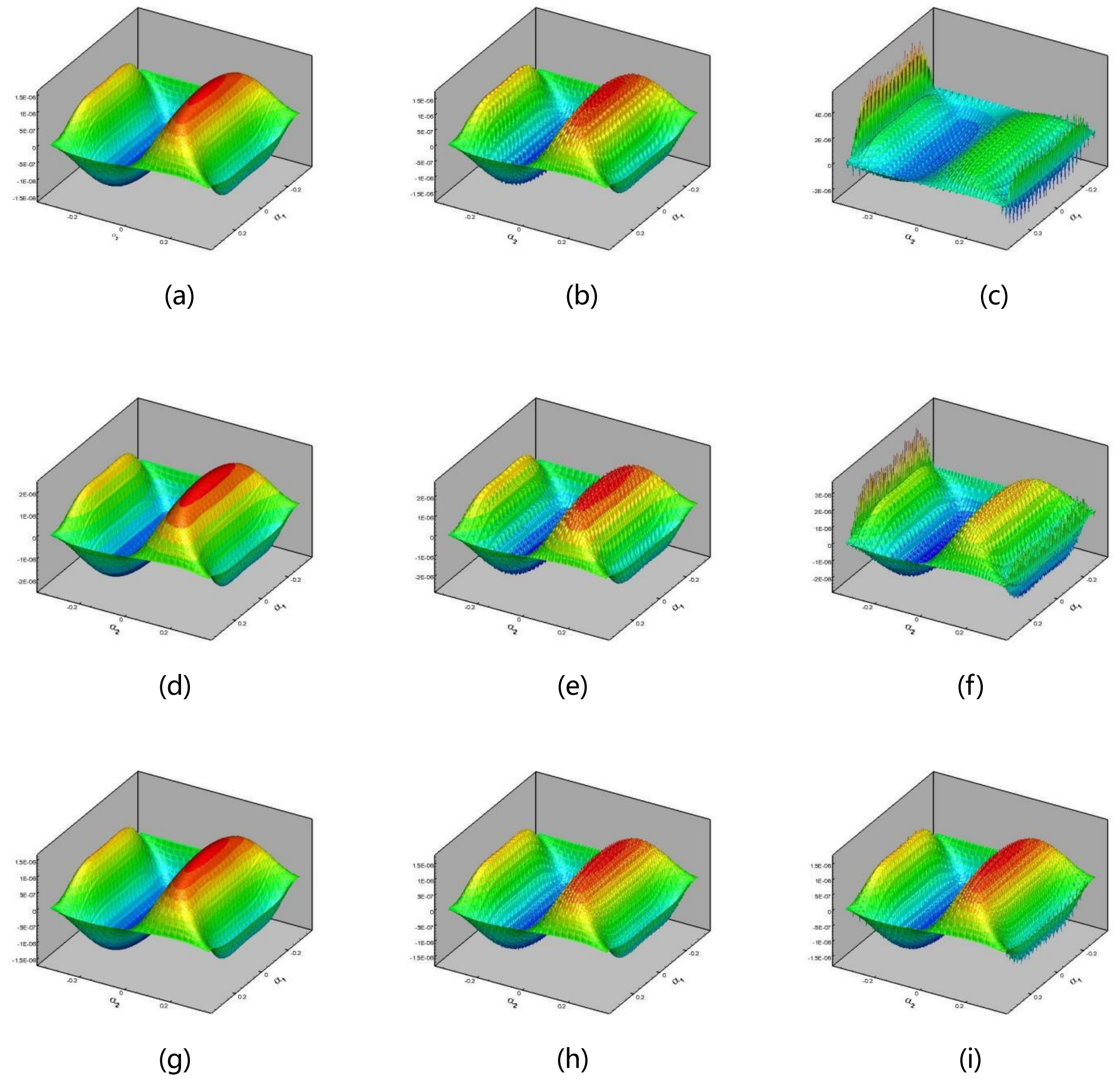

图 3 复合材料壳位移场计算结果示意图, 细观单胞 I: (a) $\boldsymbol{u}_2^{(0)}$; (b) $\boldsymbol{u}_2^{(1\varepsilon)}$; (c) $\boldsymbol{u}_2^{(2\varepsilon)}$; 细观单胞 II: (d) $\boldsymbol{u}_2^{(0)}$; (e) $\boldsymbol{u}_2^{(1\varepsilon)}$; (f) $\boldsymbol{u}_2^{(2\varepsilon)}$; 细观单胞 III: (g) $\boldsymbol{u}_2^{(0)}$; (h) $\boldsymbol{u}_2^{(1\varepsilon)}$; (i) $\boldsymbol{u}_2^{(2\varepsilon)}$.

Fig. 3. The computational results of displacement field of composite shell, mesoscopic configuration I: (a) $\boldsymbol{u}_2^{(0)}$; (b) $\boldsymbol{u}_2^{(1\varepsilon)}$; (c) $\boldsymbol{u}_2^{(2\varepsilon)}$; mesoscopic configuration II: (d) $\boldsymbol{u}_2^{(0)}$; (e) $\boldsymbol{u}_2^{(1\varepsilon)}$; (f) $\boldsymbol{u}_2^{(2\varepsilon)}$; mesoscopic configuration III: (g) $\boldsymbol{u}_2^{(0)}$; (h) $\boldsymbol{u}_2^{(1\varepsilon)}$; (i) $\boldsymbol{u}_2^{(2\varepsilon)}$.

材料板壳结构细观尺度的局部振荡行为, 因此实际工程计算中应采用高精度的二阶双尺度解. 而直接采用细网格上的有限元方法求解复合材料板壳结构的多尺度力学问题(1)时, 需要非常精细的网格才能够捕捉到复合材料板壳结构的局部振荡行为, 此时有限元计算无法保证收敛, 导致无法得到可靠的有限元参考解.

## 5.2 复合材料板壳结构强度预测

本节对由正交编织复合材料、曲线编织复合材料和四向编织复合材料组成的复合材料板和壳的屈服强度进行计算. 我们研究的复合材料板壳结构的基体杨氏模量为 410.00GPa, 泊松比为

0.18, 屈服强度为 0.14GPa. 夹杂的杨氏模量为 240.00GPa, 泊松比为 0.20, 屈服强度为 3.50MPa. 为了模拟飞行器结构所受气动力响应, 对复合材料板和壳上表面中心区域施加垂直面力如图 4 所示, 并根据已建立的复合材料板壳强度预测的高阶多尺度模型, 对复合材料板壳结构的屈服强度进行预测.

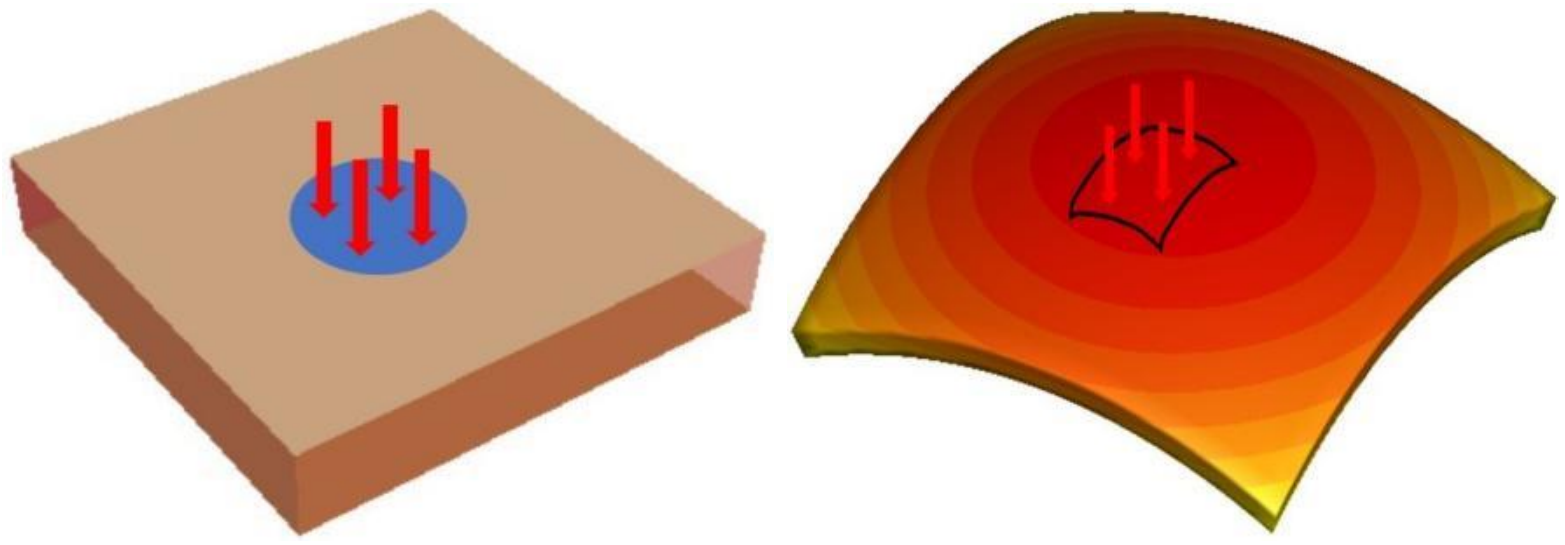

图 4 复合材料板和壳受力示意图: (a) 复合材料板结构; (b) 复合材料壳结构.

Fig. 4.Diagram of mechanical loadings on composite plate and shell: (a) macroscopic structure of composite plate; (b) macroscopic structure of composite shell.

对复合材料板和壳结构的数值模拟采用的有限元计算网格同上一节力学行为模拟时的计算网格, 数值模拟完成后, 强度预测结果如下表 2 所示.

表 2 复合材料板壳结构强度预测结果.

Table 2. The predictive results of yield strengths of composite plate and shell structures.

| 结构和单胞模型 | 均匀化方法 | 一阶双尺度方法 | 二阶双尺度方法 |
|---|---|---|---|
| 正交编织复合材料板 | 59.53MPa | 59.53MPa | 59.53MPa |
| 曲线编织复合材料板 | 48.22MPa | 48.22MPa | 48.22MPa |
| 四向编织复合材料板 | 53.92MPa | 53.92MPa | 53.92MPa |
| 正交编织复合材料壳 | 0.5736MPa | 0.5736MPa | 0.5736MPa |
| 曲线编织复合材料壳 | 0.4745MPa | 0.4745MPa | 0.4745MPa |
| 四向编织复合材料壳 | 0.5254MPa | 0.5254MPa | 0.5254MPa |

这里需要强调均匀化方法、一阶双尺度方法和二阶双尺度方法预测的复合材料板和壳结构的强度一致是因为其基体材料和夹杂材料的材料参数差异性不大, 此时复合材料板和壳结构的细观尺度振荡行为不明显. 但对于基体材料和夹杂材料材料参数差异大时, 使用高阶多尺度方法预测其屈服强度是十分必要的. 此外, 我们利用 ANSYS 软件计算了复合材料板和壳位移场的均匀化解, 与我们采用统一多尺度方法计算得到的复合材料板和壳位移场的均匀化解基本一致, 验证了我们的统一多尺度方法的可靠性, 位移计算结果如下表 3 所示.

表 3 复合材料板壳结构位移场计算结果.
Table 3. The computational results of displacement field of composite plate and shell structures.

| 结构和单胞模型 | 位移场极值(ANSYS 软件) | 位移场极值(均匀化方法) |
|---|---|---|
| 正交编织复合材料板 | -1.8155e-6m | -1.6231e-6m |
| 曲线编织复合材料板 | -1.7914e-6m | -1.6137e-6m |
| 四向编织复合材料板 | -1.8113e-6m | -1.6188e-6m |
| 正交编织复合材料壳 | -5.1033e-7m | -5.8758e-7m |
| 曲线编织复合材料壳 | -5.1310e-7m | -5.3218e-7m |
| 四向编织复合材料壳 | -5.1421e-7m | -5.3187e-7m |

# 6 结 论

本文建立了复合材料板壳结构力学行为模拟和强度预测的统一高阶多尺度方法, 建立的高阶多尺度模型将复合材料板力学行为的分析统一到复合材料壳的分析. 此外, 建立的高阶多尺度模型保证了力学平衡方程的小尺度逼近性, 可以精确地捕捉复合材料板壳结构细观尺度的振荡行为. 基于建立的复合材料板壳结构力学行为的高阶多尺度模型, 结合材料强度理论, 建立了复合材料板壳结构强度预测的高阶多尺度模型. 最后结合有限元方法和插值方法, 建立了复合材料力学行为模拟和强度预测的高效多尺度算法. 数值实验表明, 本文发展的高阶多尺度方法可以高精度、高效率模拟复合材料板壳结构的力学行为并能够有效预测复合材料板壳结构的屈服强度. 本文建立的统一的高阶多尺度计算框架可以很容易地扩展到多场耦合环境下复合材料板壳结构力学性能的计算和强度预测, 为复合材料板壳结构的结构设计和性能优化提供了坚实的基础理论和高性能算法. 在未来, 我们将利用建立的高阶多尺度计算框架对复合材料板壳结构的其他非线性物理、力学性能进行深入研究.

# 7 致谢

感谢中国科学院数学与系统科学研究院崔俊芝院士的讨论.

# Unified high-order multi-scale method for mechanical behavior simulation and strength prediction of composite plate and shell structures*

Dong Hao[1)]† Ge Bu-Feng[2)] Gao Ming-Yuan[3)]

1)(*School of Mathematics and Statistics, Xidian University, Xi' an 710071, China*)

2)(*School of Electronic Engineering, Xidian University, Xi' an 710071, China*)

3)(*School of Artificial Intelligence, Xidian University, Xi' an 710071, China*)

# 8 Abstract

The complicated mesoscopic configurations of composite plate and shell structures requires a huge amount of computational overhead for directly simulating their mechanical problems. In this paper, a unified high-order multi-scale method, which can effectively simulate the mechanical behavior and predict yield strength of composite plates and shells, is developed. Firstly, through the multiscale asymptotic analysis of multi-scale elastic equations in the orthogonal curvilinear coordinate system, a high-order multi-scale model is established, which can uniformly and effectively analyze the mechanical behavior of composite plate and shell structures. Moreover, the error estimation of the high-order multi-scale solutions is derived. Then, combining with the material strength theory, a high-order multi-scale model for the strength prediction of composite plate and shell structures is established. Next, based on the established high-order multi-scale model, a multi-scale algorithm is developed which can not only efficiently and accurately simulate the mechanical behaviors of composite plate and shell structures, but also predict their yield strength. Finally, the effectiveness of the established high-order multi-scale method is verified by extensive numerical experiments. The numerical experimental results indicate that the high-order multi-scale method can more accurately capture the meso-scale oscillatory behaviors of composite plate and shell structures. The unified high-order multi-scale method established in this paper is not only suitable for the prediction of mechanical properties of composite plate and shell structures, but also can be further extended to the prediction of multi-field coupling properties of composite plate and shell structures.

Keywords: Compositeplatesand shells, Multi-scale asymptotic analysis, Theory of strength, Highorder multi-scale method, Error analysis

PACS: 72.80.Tm, 47.11.St, 46.50.+a, 02.60.-x

* Project supported by the National Natural Science Foundation of China (Grant No. 12001414), the Fundamental Research Funds for the Central Universities (No. JB210702), the open foundation of Hubei Key Laboratory of Theory and Application of Advanced Materials Mechanics (Wuhan University of Technology) (No. WUT-TAM202104).

† Corresponding author. E-mail: donghao@mail.nwpu.edu.cntel: 13468772907